
\let\1\sp  
\catcode`@=11
\newcount\refno \refno=1
\newdimen\refskip  \refskip=\parskip
\def\eat#1{}
\def\ifundefined#1{\expandafter\ifx\csname
                        \expandafter\eat\string#1\endcsname\relax}
\def\ifempty#1{\ifx\@mp#1\@mp}
\def\doref  #1 #2#3\par{{\refno=0
  \ifundefined#2\def#2{{\tt\string#2}}\fi
     \vbox {\everyref \item {\reflistitem#2}
            {\d@more#3\more\@ut\par}\par}}\vskip\refskip }
\def\d@more #1\more#2\par
   {{#1\more}\ifx#2\@ut\else\d@more#2\par\fi}
\def\@utfirst #1,#2\@ver
   {\author#1,\ifx#2\@ut\afteraut\else\@utsecond#2\@ver\fi}
\def\@utsecond #1,#2\@ver
   {\ifx#2\@ut\andone\author#1,\afterauts\else
      ,\author#1,\@utmore#2\@ver\fi}
\def\@utmore #1,#2\@ver
   {\ifx#2\@ut\and\author#1,\afterauts\else
      ,\author#1,\@utmore#2\@ver\fi}
\def\authors#1{\@utfirst#1,\@ut\@ver}
\catcode`@=12 

\let\REF\labref
\def\cite#1{\citeform{#1}}  

\def\citeform#1{{\bf\lbrack#1\rbrack}}
\let\everyref\relax            
\let\more\relax                
\let\reflistitem\citeform

\def\Bref#1 "#2"#3\more{\authors{#1}:\ {\it #2}, #3\more}
\def\Gref#1 "#2"#3\more{\authors{#1}\ifempty{#2}\else:``#2''\fi,
                             #3\more}
\def\Jref#1 "#2"#3\more{\authors{#1}:``#2'', \Jn#3\more}
\def\inPr#1 "#2"#3\more{in: \authors{\eds#1}:``#2'', #3\more}
\def\Jn #1 @#2(#3)#4\more{{\it#1}\ {\bf#2}(#3)#4\more}
\def\author#1. #2,{#1.~#2}
\def\sameauthor#1{\leavevmode$\underline{\hbox to 25pt{}}$}
\def\and{, and}   \def\andone{ and}
\let\afteraut\relax
\let\afterauts\relax
\def\etal{\def\afteraut{, et.al.}\let\afterauts\afteraut
           \let\and,}
\def\eds{\def\afteraut{(ed.)}\def\afterauts{(eds.)}}
\overfullrule=0pt
\magnification 1200
\baselineskip=12pt
\hsize=16.5truecm \vsize=23 truecm \voffset=.4truecm
\parskip=14 pt

\def\capital#1{\expandafter\capitalula#1lalula}
\def\capitalula#1#2lalula{\uppercase{#1}#2}
\def\bgsection#1{\vskip0pt plus.1\vsize\penalty-250
  \vskip0pt plus-.1\vsize\bigskip\vskip\parskip
  \message{#1}\leftline{\bf#1}\nobreak\smallskip\noindent}
\def\ACKNOW#1\par{\ifx\REF\doref \bgsection{Acknowledgements}
                    #1\par\bgsection{References}\fi}

\def\idty{{\leavevmode{\rm 1\ifmmode\mkern -5.4mu\else\kern -.3em\fi I}}}
\def\Ibb #1{ {\rm I\ifmmode\mkern -3.6mu\else\kern -.2em\fi#1}}
\def\Ird{{\hbox{\kern2pt\vbox{\hrule height0pt depth.4pt width5.7pt
    \hbox{\kern-1pt\sevensy\char"36\kern2pt\char"36} \vskip-.2pt
    \hrule height.4pt depth0pt width6pt}}}}
\def\Irs{{\hbox{\kern2pt\vbox{\hrule height0pt depth.34pt width5pt
       \hbox{\kern-1pt\fivesy\char"36\kern1.6pt\char"36} \vskip -.1pt
       \hrule height .34 pt depth 0pt width 5.1 pt}}}}

\def\ibb #1{\leavevmode\hbox{\kern.3em\vrule
     height 1.5ex depth -.1ex width .2pt\kern-.3em\rm#1}}
\def\Nl{{\Ibb N}} \def\Cx {{\ibb C}} \def\Rl {{\Ibb R}}

\def\lessblank{\parskip=5pt \abovedisplayskip=2pt
          \belowdisplayskip=2pt }
\outer\def\iproclaim #1. {\vskip0pt plus50pt \par\noindent
     {\bf #1.\ }\begingroup \interlinepenalty=250\lessblank\sl}
\def\eproclaim{\par\endgroup\vskip0pt plus100pt\noindent}
\def\proof#1{\par\noindent {\bf Proof #1}\          
         \begingroup\lessblank\parindent=0pt}
\def\QED {\hfill\endgroup\break
     \line{\hfill{\vrule height 1.8ex width 1.8ex }\quad}
      \vskip 0pt plus100pt}


\def\Bar{\overline}
\def\abs #1{{\left\vert#1\right\vert}}
\def\bra #1>{\langle #1\rangle}
\def\bracks #1{\lbrack #1\rbrack}
\def\dim{\mathop{\rm dim}\nolimits}

\def\ket #1{\vert#1\rangle}

\def\norm #1{\left\Vert #1\right\Vert}

\def\rank{\mathop{\rm rank}\nolimits}
\def\set #1{\left\lbrace#1\right\rbrace}
\def\stt{\,\vrule\ }
\def\Set#1#2{#1\lbrace#2#1\rbrace}  

\def\rstr{\hbox{$\vert\mkern-4.8mu\hbox{\rm\`{}}\mkern-3mu$}}

\def\phi{\varphi}
\def\epsilon{\varepsilon}
\def\3{\ss}
\def\Re{\mathchar"023C\mkern-2mu e}

\def\em{\it} 
\def\SnU#1{S_\nu U(#1)} 
\def\B{{\cal B}} 
\def\C{{\cal C}} 
\def\N{{\cal N}} 
\def\H{{\cal H}}         
\def\R{{\cal R}}         
\def\cE{{\cal E}}            
\def\EH#1{\cE_{\H}(#1)}
\def\EHq{\EH q}
\def\Cuntzal{{\cal O}}
\def\OH#1{\Cuntzal_{\H}(#1)}
\def\OHq{\OH q}
\def\Exp{{\rm Exp}}       
\def\absq{\abs q}
\def\M{{\cal M}}    
\def\as{\a\1\dagger}
\def\a{a}
\def\vs{v\1\dagger}
\def\isom{\eta} 
\def\isOm{\widehat\eta} 
\def\thetah{\widehat\theta}  
\def\St{{\rm S}}             
\def\rep{representation}
\def\irrep{irreducible \rep}
\def\pco{peripheral coherent}
\def\refdefs #1,#2\last{\global\edef#1{\number\refno}
                        \global\advance\refno by1
                         \ifx#2\next\else\refdefs#2\last\fi}
\refdefs \Askey,\Baez,\Bergmann,\Bieden,\Bourbaki,\BoSpeia,
    \JOa,\BraRo,\Cuntz,\Dix,\Dunford,\Dykema,\Fivel,\Fuglede,
    \Green,\Greenloc,\Rieckers,\QCR,\GPOTS,\WICK,\Katriel,\next\last
\refdefs \Klauder,\KuSpei,\Burkhard,\Lawson,\LiSheng,\Farlane,
    \Nica,\Foias,\Rudin ,\Sakai,\Shale,\Speia,\QTD,
    \Worono,\Zagier,\next\last
\line{}
\vskip 2.0cm
\begingroup\nopagenumbers
\font\BF=cmbx10 scaled \magstep 3
\newcount\fnotes \fnotes=1
\def\Footnote#1{\footnote{$\1{\number\fnotes}$}%
{#1}\advance\fnotes by1}

{\BF \baselineskip=25pt
\centerline{     Coherent States of the }
\centerline{     q--Canonical Commutation Relations }
}
\vskip 1.0cm

\centerline{
{\bf P.E.T. J\o rgensen
\Footnote
  {{\sl Dept. of Mathematics,
        University of Iowa,
        Iowa City, IA 52242, USA}
}$\1,$\Footnote
    {{\sl Supported in part by the NSF(USA), and NATO}},
\bf and R.F. Werner
\Footnote{{\sl FB Physik, Universit\"at Osnabr\"uck,
            Postfach 4469, \hfill\break\vrule width0pt\qquad
          D-4500 Osnabr\"uck, Germany.
}}$\1,$\Footnote
{{\sl Electronic mail:\quad \tt reinwer@dosuni1.rz.Uni-Osnabrueck.DE}}
}}
\vskip 1.0cm
\vfil

{\baselineskip=12pt
\midinsert\narrower\narrower\noindent
{\bf Abstract.}
For the $q$-deformed canonical commutation relations\break
$a(f)a\1\dagger(g)
        =(1-q)\,\langle f,g\rangle\idty+q\,a\1\dagger(g)a(f)$
for $f,g$ in some Hilbert space ${\cal H}$
we consider representations generated from a vector $\Omega$
satisfying $a(f)\Omega=\langle f,\phi\rangle\Omega$,
where $\phi\in{\cal H}$.
We show that such a representation exists if and only if
$\Vert\phi\Vert\leq1$. Moreover, for $\Vert\phi\Vert<1$ these
representations are unitarily equivalent to the Fock representation
(obtained for $\phi=0$). On the other hand representations obtained
for different unit vectors $\phi$ are disjoint.
We show that the universal C*-algebra for the relations has a
largest proper, closed, two-sided ideal. The quotient by this ideal
is a natural $q$-analogue of the Cuntz algebra (obtained for $q=0$).
We discuss the Conjecture that, for $d<\infty$, this analogue
should, in fact, be equal to the Cuntz algebra itself.
In the limiting cases $q=\pm1$ we determine all irreducible
representations of the relations, and characterize those which can
be obtained via coherent states.
\endinsert
}
\vskip 1truecm


\noindent
{\bf AMS subject classification: }
46L89, 58B30, 81R10, 81R30

\vskip 1cm
\vfil\eject\endgroup
\beginsection 1. Introduction

In this paper we study some new aspects of a set of commutation
relations, depending on a parameter $q\in(-1,1)$ studied by various
authors on quite different motivations. Greenberg \cite\Green\
introduced these relations as an interpolation between Bose ($q=1$)
and Fermi ($q=-1$) statistics. He was particularly interested in the
observable consequences of a hypothetical small deviation from the
Pauli principle. However, due to problems with field theoretical
localizability \cite\Greenloc\ and thermodynamic stability
\cite\QTD, a naive particle interpretation of systems satisfying
these relations is problematic.
Speicher \cite\Speia\ introduced these relations as a
new kind of quantum ``noise'', which could be used as a driving force
in a quantum stochastic differential equation \cite\KuSpei. From the
point of view of C*-algebra theory the relations became interesting as
an example of a C*-algebra defined in terms of generators and
relations. In this context it was observed that the relations reduce
for $q=0$ to those studied by Cuntz \cite\Cuntz.

The special case of a single generator, the so-called $q$-oscillator,
was introduced by Biedenharn \cite\Bieden\ and Macfarlane
\cite\Farlane\ as a means of constructing \rep s of quantum
groups. In fact, the $q$-oscillator also appears as a subalgebra of
the quantum group $\SnU2$ \cite\Worono. The $q$-oscillator can be
studied in full detail by representing the generator as a weighted
unilateral shift (in mathematical terminology) or as a Bose creation
operator multiplied with a suitable function of the number operator
(in physical terminology). This has been noted in a large number of
papers. We will use this representation in the present paper to obtain
information about the non-trivial case of several generators.

In this case most early work \cite{\Green,\BoSpeia,\Fivel} focussed on
showing that the scalar product in the $q$-analogue of
Fock space is positive definite.
On the other hand, from the C*-algebraic point of view the most
immediate and natural problem arising from the relations was to
characterize the norm-closed operator algebra generated by any
realization of the relations by bounded operators on a Hilbert
space. Here the case $q=0$ served as a model: for $q=0$ this algebra
must be either isomorphic to the one obtained in the Fock \rep,
called the Cuntz-Toeplitz algebra, or a quotient of the
Cuntz-Toeplitz algebra by its unique two-sided ideal (isomorphic to
the compact operators), known as the Cuntz algebra. For $q\neq0$ the
first important step was made in \cite\QCR, where we showed that for
$\abs q<\sqrt2-1\approx .41$ the same results hold. In particular,
the C*-algebras generated with $q$ in this range are exactly the
same as for $q=0$. The condition $\abs q<\sqrt2-1$ is certainly not
optimal, and all results known to us are compatible with the
conjecture (which we will refer to as ``{\it Conjecture C\/}'', see
Section 4) that the results of \cite\QCR\ hold for all $\abs q<1$.
However, no decisive progress towards proving this conjecture has
been made since \cite\QCR. Based on an improved understanding of the
Fock \rep\ \cite\Zagier, Dykema and Nica \cite\Dykema\ managed to
extend the interval for $q$ slightly, but only for the algebra
generated in the Fock \rep. More importantly, they established, for
the Fock \rep\ only, the existence of the homomorphism between the
algebras for $q=0$ and for general $-1<q<1$, which according to
Conjecture C should be an isomorphism. We will briefly describe and
apply their results in Section 4.

The main aim of this paper is to study the $q$-analogue of a structure
which is well-known in the limiting cases $q=0,\pm1$, namely the
generalization of the Fock state to the so-called coherent states.
In the case of a single relation such states appear in \cite\LiSheng,
although, due to a different choice of generators, their work makes
sense only in the Fock \rep, and gives states different from ours.
We will determine all coherent states, and discuss under what
circumstances they generate the same \rep, or are mutually singular.
Using coherent states, we show that the universal C*-algebra
generated by the relations has a unique largest closed two-sided
ideal. (If Conjecture C holds this ideal is also the only proper
ideal, and isomorphic to the compact operators). The quotient of the
algebra by this ideal is then simple, and the natural analogoue of
the Cuntz algebra for $q\neq0$. Finally, we consider the limiting
case $q=-1$, and compute all \irrep s of the relations with Clifford
algebra methods. It turns out that in this degenerate case the
coherent states exhaust only a small subclass of \irrep s.

We emphasize that when we talk about \rep s in the sequel we always
mean *-\rep s of some involutive algebra by {\it bounded\/}
operators on a Hilbert space. Thus even if the relations may have
interesting unbounded realizations we do not consider them.

\beginsection 2. $q$-relations and coherent states

The following Proposition introduces the ``$q$-relations'' which are
the object of our study.

\iproclaim 1 Proposition.
\def\ast{\widetilde a\1\dagger}%
Let $\H$ be a Hilbert space, and let $q\in\Rl$, $\absq<1$.
Then there is a C*-algebra $\EHq$ generated by elements
$\as(f)$ for $f\in\H$, such that $f\mapsto\as(f)$ is linear, and
$$  \a(f)\as(g)=(1-q)\,\bra f,g>\idty+q\,\as(g)\a(f)
\quad,\eqno(1)$$
where $\a(f):=\as(f)\1*$.
For $q=1$, and orthogonal unit vectors $e_1,\ldots,e_n\in\H$ the bound
$$ \sum_{i=1}\1n \as(e_i)\a(e_i)\leq\idty
\eqno(2)$$
holds.
\hfill\break
Moreover, $\EHq$ is uniquely determined by the following universal
property: whenever $\widetilde{\cE}$ is a C*-algebra containing
elements $\ast(f)$ satisfying the above conditions, there
is a unique unital homomorphism $\phi:\EHq\to\widetilde{\cE}$ such that
$\phi(\as(f))=\ast(f)$.
\eproclaim

The proof of this result is given in \cite\QCR. Note that in
comparison with \cite{\QCR,\BoSpeia} we have changed the
normalization of the operators $\as(f)$. This modification makes no
essential difference for $\absq<1$. However, it removes the
singularity of the relations for $q\to1$, and simplifies all
algebraic expressions.
Moreover, it was shown in \cite\Nica\ that with this normalization
the algebras $\cE_\Cx(q)$ form a continuous field of C*-algebras
\cite\Dix.
We may consider the relations (1) for all
$q\in\Rl\cup\set\infty$, where for $q=\infty$ we set
$\as(g)\a(f)=\bra f,g>\idty$. The study of the case $\absq\geq1$ is
then reduced to the case $\absq\leq1$  by the symmetry
$$\eqalign{
      q&\mapsto q\1{-1}\cr
 \as(f)&\mapsto \a(\Bar f)
}\eqno(3)$$
for some antiunitary operator $f\mapsto\Bar f$.

The crucial feature of the relations (1) is that they allow us to
order any polynomial in the generators in such a way that in every
monomial all operators $\a(f)$ are to the right of every $\as(f)$. This
normal ordered, or ``Wick ordered'' form of a polynomial is unique
\cite{\Baez,\Bergmann,\WICK},
hence we can define a linear functional $\omega$ on the polynomial
algebra over the relations by choosing an arbitrary multilinear
expression for $\omega(\as(f_1)\cdots \as(f_n)\a(g_1)\cdots
\a(g_m)\bigr)$. Since such monomials generate $\EHq$ this is also a
way to parametrize all states on the C*-algebra $\EHq$. The
following Theorem introduces the coherent states on $\EHq$ using
such a parametrization.

\iproclaim 2 Theorem.
Let $\absq\leq1$, and $\phi\in\H$ with $\norm{\phi}\leq1$.
Then there is a unique state $\omega_\phi$ on $\EHq$ such that
$$ \omega_\phi(\as(f)X)=\bra\phi,f>\, \omega(X)
\eqno(4)$$
for all $f\in\H$, and all $X\in\EHq$.
The state $\omega_\phi$ is pure.
For $\norm{\phi}>1$, there is no state satisfying (4).
\eproclaim

We will call $\omega_\phi$ the {\em coherent state} associated with
$\phi$. This terminology originated in quantum optics, where these
states are used to describe states of the electromagnetic
radiation field \cite{\Klauder,\Rieckers}. The special state
$\omega_0$ is called the {\em Fock state}. If $\norm{\phi}=1$, we will
call $\omega_\phi$ a {\em \pco} state. For any $\phi$ we will denote
by $\pi_\phi$ the GNS-\rep\ associated with $\omega_\phi$, and call it
the {\em coherent \rep} associated with $\phi$. For the special case
$q=0$, coherent states in this sense have been studied in \cite\JOa.

The proof of Theorem 2 is based on an analysis of the case of a
single relation. We summarize the relevant results in the following
Lemma. The assumption that $a$ is bounded is essential for this
result, i.e.\ there are also unbounded operators satisfying the
relation, and the conclusion fails for these.

\iproclaim 3 Lemma.
Let $\absq<1$, and let $a\equiv \a(e_1)$ with $\norm{e_1}=1$ be a
bounded operator on a Hilbert space $\R$ satisfying the relation
$aa\1*=(1-q)\idty+qa\1*a$.
\item{(1)}
Then $a$ is reduced by a unique decomposition
$\R\cong(\R_0\otimes\R')\oplus\R_1$, such that
\itemitem{(a)}
if $\R_1\neq\set0$, $a\rstr\R_1$ is unitary.
\itemitem{(b)}
if $\R'\neq\set0$, $a\rstr\R_0\otimes\R'$ acts as
$a=a_0\otimes\idty$, where $a_0$ is given explicitly as the weighted
shift
$$ a_0\1*\ket{n}= (1-q\1{n+1})\1{1/2}\, \ket{n+1}
\quad,\eqno(5)$$
where $\ket{n}$ for $n=0,1,\ldots$ is an orthonormal basis of
$\R_0$.

\item{(2)}
$$      \norm{a}=\cases{  1        & for $q>0$, or $a$ unitary \cr
                        \sqrt{1-q} & for $q<0$, and $a$ not unitary\cr}
\quad.\eqno(6)$$
\item{(3)}
There are functions $\beta_+(q)<\infty$ and $\beta_-(q)>0$ such that
$$    \beta_-(q)\idty \leq a\1n(a\1*)\1n \leq  \beta_+(q)\idty
\quad,\eqno(7)$$
uniformly for $n\in\Nl$. In particular, the spectral radius of $a$
is equal to $1$.
\item{(4)}
Let $a\1*\xi=\lambda\xi$ for $\xi\neq0$. Then
$\xi\in\R_1$, $\abs\lambda=1$, and $a\xi=\Bar\lambda\xi$.
\eproclaim

\proof:
For (1) see \cite\QCR; for (2) see \cite{\BoSpeia,\QCR}.

(3) For the unitary part $a\rstr\R_1$ we only need
$\beta_-(q)\leq1\leq\beta_+(q)$, which will be true
for the $\beta_\pm$ constructed below. Hence it suffices to take
$a=a_0$. Then
$$ a_0\1n(a_0\1*)\1n \ket{k}
      =\lambda_{k+1}\cdots\lambda_{k+n}\ket{k}
\quad,$$
where $\lambda_k=(1-q\1k)$. We will take $\beta_\pm$ as the supremum
(resp. infimum) over all products $\prod_{k\in M} \lambda_k$  for
$M\subset\Nl$. Explicitly,
$$\eqalign{
    \beta_+(q)&=\cases{1 & $q\geq0$ \cr
                       \prod_{k=1}\1\infty \bigl(1-q\1{2k+1}\bigr)
                         & $q\leq0$}     \cr
    \beta_-(q)&=\cases{\prod_{k=1}\1\infty \bigl(1-q\1{k}\bigr)
                         & $q\geq0$   \cr
                       \prod_{k=1}\1\infty \bigl(1-q\1{2k}\bigr)
                         & $q\leq0$}
\cr}\eqno(8)$$
Since these products (related to Theta functions, and
to ``$q$-factorials'' \cite\Askey) are absolutely convergent,
$\beta_\pm(q)$ is finite and non-zero for all $q,\ \absq<1$.
For computing the spectral radius we let $n\to\infty$ in the
inequality
$$       \beta_-(q)\1{1/2n} \leq \norm{a\1n}\1{1/n}
                           \leq \beta_+(q)\1{1/2n}
\quad.$$

(4)
Given the decomposition it suffices to show that
$a_0\1*\xi=\lambda\xi$ implies $\xi=0$. This follows immediately from
the weighted shift structure (5) of $\a_0\1*$, by solving the recursion
for the coefficients $\xi_n$ in $\xi=\sum_n\xi_n\ket{n}$.
\QED

Consider the GNS-\rep\ $\pi_\phi$ associated with the
coherent state $\omega_\phi$. This has a cyclic vector $\Omega_\phi$,
which is a joint eigenvector of the generators, i.e.\
$$ \a(f)\Omega_\phi= \bra f,\phi>\,\Omega_\phi
\quad.\eqno(9)$$
Conversely, any unit vector satisfying (9) will give the coherent
state via $\omega_\phi(X)$ $=\bra\Omega_\phi,X\Omega_\phi>$.
Therefore, in order to show that $\omega_\phi$ is positive, it is
sufficient to exhibit such a vector in a \rep\ which is
known to be positive. Now the Fock \rep\ $\pi_0$ has been
proven to be positive \cite{\BoSpeia,\Fivel,\Zagier,\WICK}.
Hence it suffices to find such vectors in the Fock
\rep. The basic construction for such vectors can be
carried out in the case of a single generator. For Boson commutation
relations the operator transforming the vacuum into a coherent state
is well-known to be $\exp\bigl(za\1*\bigr)$. For the $q$-relations a
similar role is played by the ``$q$-exponential'' function $\Exp_q$,
defined by the functional equation \cite\Katriel
$$ {\bf D}_q\, \Exp_q(z)
   \equiv {\Exp_q(z)-\Exp_q(q\,z)  \over z-qz}
   = \Exp_q(z)
\quad.\eqno(10)$$
The $q$-exponential satisfies no simple addition formula, and
therefore the operator connecting different coherent state can only
be expressed as a quotient of two such exponentials. Rather than
defining first the $q$-exponential, and then studying its
invertibility, we define, in the following Lemma, all these quotients
at the same time. The connection with the $q$-exponential is
$V_{\alpha0}(z)=E_q\bigl(\alpha z/(q-1)\bigr)$.

\iproclaim 4 Lemma.
\item{(1)}
Let $\absq<1$, and $\alpha,\beta\in\Rl$. Then the functional equation
$$ V_{\alpha\beta}(qz)
   = {1-\alpha z\over 1-\beta z}V_{\alpha\beta}(z)
\qquad;\quad V_{\alpha\beta}(0)=1
\eqno(11)$$
has a unique analytic solution near $z=0$, which is analytic for
$\abs{\alpha z}<1$.
For $\abs{\alpha z}<1$, and $\abs{\beta z}<1$, and $\gamma\in\Rl$ we
have $V_{\alpha\beta}V_{\beta\gamma}=V_{\alpha\gamma}$.
\item{(2)}
Let $a$ be a bounded operator on a Hilbert space $\R$ with
$aa\1*=(1-q)\idty+qa\1*a$. Then, for $\Omega_\beta\in\R$, and
$\abs\alpha<1$ we have the implication
$$ \bigl(a-\beta\bigr) \Omega_\beta=0
\qquad\Rightarrow\quad
   \bigl(a-\alpha\bigr)V_{\alpha\beta}(a\1*)\Omega_\beta =0
\quad,\eqno(12)$$
where the function $V_{\alpha\beta}$ is evaluated on $a\1*$ in the
analytic functional calculus.
\eproclaim

\proof:
Let $V_{\alpha\beta}(z)=\sum_kc_kz\1k$. Then equation (12) together
with the iterated relation
$$ a(a\1*)\1k= q\1k (a\1*)\1k a  + (1-q\1k)(a\1*)\1{k-1}
\eqno(13)$$
gives a functional equation for the coefficients $c_k$:
$$ c_{k+1}={\alpha-q\1k\beta  \over 1-q\1{k+1}} c_k
\qquad;\quad c_0=1.
\quad.\eqno(14)$$
By an elementary computation this is the same recursion which holds
for the coefficients of $V_{\alpha\beta}$ defined through the
functional equation. By standard theorems on power series its radius
of convergence is $\abs\alpha\1{-1}$.
The chain relation $V_{\alpha\beta}V_{\beta\gamma}=V_{\alpha\gamma}$
follows directly from the functional equation.
\QED

\proof{of Theorem 2:}
Let $\omega$ be a state satisfying equation (4). Then we can compute
it on any polynomial in the generators by Wick ordering the
polynomial, and then applying successively equation (4) and its
adjoint $\omega(X\a(g))=\bra g,\phi>\omega(X)$. Since polynomials are
dense in $\EHq$, $\omega=\omega_\phi$ is uniquely determined.
It is also clear that $\omega_\phi$ must be a pure state, since it is
the only state on which the positive elements
$(\as(f)-\bra\phi,f>\idty)(a(f)-\bra f,\phi>\idty)$ have zero
expectation for all $f\in\H$.

If there is a state $\omega_\phi$ we have,
for $\norm{f}=1$:
$ \abs{\bra\phi,f>}\12
    =\omega_\phi\bigl(\as(f)\1N\a(f)\1N\bigr)\1{1/N}
    \leq1$,
since the spectral radius of $\a(f)$ is $1$ by Lemma 3(2). With
$f=\phi/\norm{\phi}$ this implies $\norm{\phi}\leq1$.

It remains to be shown that $\omega_\phi(X\1*X)$ is positive for
$\norm{\phi}\leq1$, and $X\in\EHq$. Since polynomials in the
generators are norm dense in $\EHq$ by the universal property, it
suffices to show this for polynomials $X$. For such $X$,
$\omega_\phi(X\1*X)$ is obviously a continuous function of $\phi$.
Hence it suffices to show positivity for $\norm{\phi}<1$.

Let $\norm{\phi}<1$. We know from \cite\BoSpeia\ that
$\omega_0$, the Fock state, is positive.
By Lemma 3(2), $\as(\phi)$ has spectral radius $<1$. Hence we can
apply $V_{10}$ from Lemma 4 to $\as(\phi)$ in the analytic functional
calculus. Let
$V\equiv V_{10}(\as(\phi))
   =V_{\norm{\phi},0}(\as(\phi/\norm{\phi}))$.
Then since $V_{0,\norm{\phi}}(\as(\phi/\norm{\phi}))=V\1{-1}$ we have
that $\Omega_\phi=V\Omega_0$ is non-zero. By Lemma 4 we have
$\a(\phi)\Omega_\phi=\bra\phi,\phi>\Omega_\phi$. On the other hand,
when $\psi\perp\phi$, we have
$\a(\psi)\as(\phi)\1n\Omega_0
   =q\1n\as(\phi)\1n\a(\psi)\Omega_0
   =0$.
With the series expansion for $V$ we find $\a(\psi)\Omega_\phi=0$.
Combining this with the result for $\phi=\psi$ we get
$\a(\psi)\Omega_\phi=\bra\psi,\phi>\Omega_\phi$, and
$\omega_\phi(X\1*X)=\bra\Omega_\phi,X\1*X\Omega_\phi>
                  =\norm{X\Omega_\phi}\12
                  \geq0$.
\QED

The proof gives more information than just the positivity of
$\omega_\phi$: by composing the operators $V_{10}(\as(\phi))$ and
$V_{10}(\as(\psi))\1{-1}$ we get the following consequence:

\iproclaim 5 Corollary.
For $\norm{\phi},\norm{\psi}<1$, the states $\omega_\phi$ and
$\omega_\psi$ are connected by an invertible element
$v_{\phi\psi}\in\EHq$ via
$\omega_\phi(X)=\omega_\psi(v_{\phi\psi}\1*Xv_{\phi\psi})$.
\eproclaim

The operators $v_{\phi\psi}$ are Radon-Nikodym derivatives in the
sense of \cite\Sakai. Since they are defined by norm convergent
series, they end up in the C*-algebra $\EHq$, and not merely in some
bigger von Neumann algebra.

We close this section with a brief discussion of the coherent states
for certain variations of the $q$-relations found in the literature.
Most of the literature is concerned with the Fock \rep\ of the
relations with a single generator, and the relations are frequently
written in a form explicitly involving the number operator $N$ of the
Fock \rep. This operator is defined by
$\exp(itN)\as(f)\exp(-itN)=\as\bigl(\exp(it)f\bigr)$, and
$N\Omega_0=0$, where $\Omega_0$ is the Fock vacuum. We will continue
to denote by $\as(f)$ the generators with the conventions fixed in
Proposition 1. Then the generators found elsewhere are $b\1\dagger(f)$
with
$$\eqalign{ b\1\dagger(f)&=\beta\, q\1{\alpha N}\,\as(f)
                         =\beta\, \as(f) q\1{\alpha (N+1)} \cr
        b(g)b\1\dagger(f)&=\abs\beta\12 (1-q)\ \bra f,g>\
                           q\1{2\alpha(N+1)} +
                           q\1{2\alpha+1}\, b\1\dagger(f)b(g)
\quad.}\eqno(15)$$
The normalization used in this paper agrees with \cite\Nica,
implicitly with \cite\Worono, and with one of the versions
introduced by \cite\Farlane\ (written with a
different parameter $\widetilde q=q\1{-1/2}$).
In most of the papers in the bibliography
we have the convention $\beta=(1-q)\1{-1/2}, \alpha=0$.
The existence of a vector $\Psi\neq0$ in Fock space with
$$b(f)\Psi=\bra f,\phi>\Psi
\eqno(16)$$
is then equivalent to $\norm{\phi}\leq(1-q)\1{-1/2}$, and the joint
eigenvectors of the $b(f)$ are precisely those of the $\a(f)$.

On the other hand, when $\alpha>0$, the series for $\Psi$ satisfying
(16) diverges for all $f\neq0$, and no joint eigenvectors can be found.
The interesting cases are for $\alpha<0$. The coefficients of the
power series then decrease more rapidly, and the series defines an
entire function. Hence no constraint is placed on $\norm{\phi}$, and
the notion of \pco\ states makes no sense. This is related to the
fact that the relations then explicitly involve the operator $N$,
and hence make sense only in the Fock \rep. The relations appeared
for the first time (for a single generator, a so called
$q$-oscillator) in \cite\Bieden\ with $\alpha=-1/4$,
$\abs\beta\12=q\1{1/2}(1-q)\1{-1/2}$, and in \cite\Farlane\ with the
same constants, but using $\widetilde q= q\1{-1/2}$. Of potential
interest is also the case $\alpha=-1/4$, $\abs\beta\12=q(1-q)\1{-1}$
in which the relations can be expressed by an ordinary commutator,
i.e.\ $\bracks{\a(f),\as(g)}=\bra f,g>q\1{-N}$.

\beginsection 3. \capital\pco\ states

A remarkable fact about the \pco\ states, i.e.\ the coherent states
$\omega_\phi$ with $\norm{\phi}=1$, is the following: if
$\H'\subset\H$, there is a canonical embedding
$\cE_{\H'}(q)\hookrightarrow\EHq$. With respect to this embedding a
\pco\ state on $\cE_{\H'}(q)$ has a {\em unique} extension to
$\EHq$, which is also a \pco\ state. This follows readily from the
first item of the following Proposition.

\iproclaim 7 Proposition.
Let $\phi\in\H$, with $\norm{\phi}=1$. Then
\item{(1)}
$\omega_\phi$ is the uniquely characterized by the condition
$\omega_\phi\bigl((\as(\phi)-\idty)(\a(\phi)-\idty)\bigr)=0$.
\item{(2)}
For $\dim\H>1$ the kernel of the GNS-\rep\ $\pi_\phi$
contains every closed two-sided ideal of $\EHq$.
\eproclaim

\proof:
We set $a=\a(\phi)$, for short.

(1)
Let $\Omega$ denote the GNS-vector of a state $\omega$ with
$\omega\bigl((a\1*-\idty)(a-\idty)\bigr)=0$.
Then we have $a\Omega=\Omega$, and since on the subspace
generated by the $(a\1*)\1n\Omega$, $a$ is unitary, we also have
$a\1*\Omega=\Omega$. Then for any vector $\psi\in\H$, we
get
$$\eqalign{
  \a(\psi)\Omega &= \a(\psi)(a\1*)\1n\Omega  \cr
                  &= (1-q\1n)\bra\psi,\phi> (a\1*)\1{n-1}\Omega
                      +q\1n (a\1*)\1n \a(\psi)\Omega
\quad.\cr}$$
Since $\norm{(a\1*)\1n}$ is uniformly bounded we can take the limit
$n\to\infty$ on the right hand side, and obtain
$\a(\psi)\Omega=\bra\psi,\phi>\Omega$. Hence $\Omega$ implements
$\omega_\phi$.

(2)
\def\ideal{{\cal J}}%
\def\quot{\widetilde{\cal E}}%
Let $\ideal\subset\EHq$  be a closed two-sided ideal, and consider
the algebra $\quot=\EHq/\ideal$ with quotient mapping
$\eta:\EHq\to\quot$.
Since $\dim\H>1$ we know from Proposition 4 in \cite\QCR\ that
$\eta(a)\in\quot$ cannot be unitary, and consequently that the spectrum
of $\eta(a)$ contains the spectrum of $a_0$, the generator in the Fock
\rep. This is the unit disk, and hence the spectrum of $\eta(a)$
includes $1$. It follows (by compactness of the state space of a
C*-algebra) that there is a \rep\ $\widetilde\pi:\quot\to\B(\R)$ in which
$1$ is an eigenvalue of $\widetilde\pi(\eta(a))$. But then by part
(1) we have
$$ \bra\xi,\pi(\eta(X))\xi>=\omega_\phi(X)
\quad,\eqno(*)$$
where $\xi$ is the corresponding normalized eigenvector. The kernel
of $\pi_\phi$ is the set of $Y\in\EHq$ such that
$\omega_\phi(X\1*YZ)=0$ for all $X,Z\in\EHq$. By equation $(*)$ it
is now plain that
$\ideal=\ker\eta
       \subset\ker\pi\circ\eta
       \subset\ker\pi_\phi$.
\QED

The second part of this Proposition suggests the following
terminology:

\iproclaim 8 Definition.
Let $\H$ be a Hilbert space with $\dim\H>1$. Then the {\bf
$q$-Cuntz} algebra $\OHq$ over $\H$ is the quotient of $\EHq$ by its
unique largest ideal. Equivalently, $\OHq=\pi_\phi\bigl(\EHq\bigr)$
for any \pco\ \rep.
\eproclaim

Of course, for $q=0$ the $q$-Cuntz algebra is just the usual Cuntz
algebra $\Cuntzal_{\dim\H}$. For $\dim\H<\infty$ Conjecture C says
that $\OHq\cong\OH0$, and this is proven \cite\QCR\ for
$\absq<\sqrt2-1$. We will further extend this interval in Section 4,
using the results of \cite\Dykema. When $\dim\H=\infty$, one can
show that the Fock \rep\ of $\EHq$ is simple \cite\Burkhard. Hence
in that case, $\OHq$ is isomorphic to the Fock \rep\ of $\EHq$.

{}From Corollary 6 we know that the non-\pco\ \rep s are all equivalent.
For the \pco\ \rep s we know that the C*-algebras $\pi_\phi(\EHq)$ are
all equal. However, the von Neumann algebras $\pi_\phi(\EHq)''$ are
not: in the following Proposition we show that all \pco\ \rep s are
disjoint.

\iproclaim 9 Proposition.
Let $\phi,\psi\in\H$, with $\norm{\phi}=\norm{\psi}=1$, and let
$\pi:\EHq\to\B(\R)$ be any \rep.
\item{(1)}
The strong operator limit
$$  P(\phi)=
        \lim_{n\to\infty} \ {1\over n}\sum_{k=1}\1n
            \pi\bigl(\as(\phi)\bigr)\1k
\eqno(15)$$
exists, and is a self-adjoint projection.
\item{(2)}
For $\chi\in\H$:\quad
$\pi\bigl(\a(\chi)\bigr)P(\phi)
      =\bra\chi,\phi>P(\phi)$.
\item{(3)}
Let $P(0)$ denote the orthogonal projection onto the space
$$ \N=\set{\Omega\in\R \stt
          \forall \phi\in\H:\ \pi\bigl(\a(\phi)\bigr)\Omega=0 }$$
of Fock vectors. Then, for $\phi,\psi$ unit vectors in $\H$, or zero,
and for $X\in\EHq$:
$$   P(\phi)\, \pi(X)\, P(\psi) = \cases{
        \omega_\phi(X) P(\phi) & $\phi=\psi$\cr
                        0      & $\phi\neq\psi$.\cr}
$$
\eproclaim

\proof:
In the proof we will suppress the \rep\ $\pi$ for
notational convenience.
The existence of the limit (1) follows from the Mean Ergodic Theorem
(e.g.\ Corollary VIII,5.4 in \cite\Dunford, and the fact that the
powers $\as(\phi)\1k$ are uniformly norm bounded by Lemma 3.(3). Let
$\chi\in\H$. Then
$$\eqalign{
  \a(\chi){1\over n}\sum_{k=1}\1n \as(\phi)\1n
     &={1\over n}\sum_{k=1}\1n \Set\Big{
              (1-q\1k)\bra\chi,\phi>\as(\phi)\1{k-1}
              +q\1k\as(\phi)\1k \a(\chi)}    \cr
     &=\bra\chi,\phi>{1\over n}\sum_{k=1}\1n \as(\phi)\1n
        +\hbox{Rest  }                       \cr
\noalign{with the estimate}
\norm{\hbox{Rest}}
     &\leq {1\over n}\norm{\idty-\as(\phi)\1n}
         + {1\over n}\,{1\over1-\absq} \beta_+(q)
\quad,\cr}$$
and $\beta_+(q)$ from equation (8).
Taking the strong limit $n\to\infty$ we find (2). In particular, we
have $\a(\phi)P(\phi)=P(\phi)$, which implies
$P(\phi)\1*P(\phi)={\rm weak-}\lim(1/n)\sum_{k=1}\1n\a(\phi)\1kP(\phi)
                 =P(\phi)$,
and hence that $P(\phi)$ is an orthogonal projection.

To prove (3), let $X$ be a polynomial in the
generators, which we may assume to be Wick ordered. Then after
finitely many applications of (2) we find that $P(\phi)X P(\psi)$ is
equal to some factors times $P(\phi)P(\psi)$. If $\phi=\psi$(possibly
$\phi=\psi=0$), the factors add up to $\omega_\phi(X)$, and
the
result follows because $P(\phi)$ is a projection.

It remains to show that $P(\phi)P(\psi)=0$, when $\phi\neq\psi$, and
$\phi\neq0$. Since $P(\phi)$ is also the weak limit of
$(1/n)\sum_{k=1}\1n\a(\phi)\1k$ this follows from (2) and the observation that
$\lim_{n\to\infty}(1/n)\sum_{k=1}\1n\bra\phi,\psi>\1k$ vanishes,
unless $\phi=\psi$.

\QED

We can use the universal \rep\ for $\pi$. Then the projections
$P(\phi)$ are interpreted as projections in the universal enveloping
algebra $\EHq\1{**}$. By (4) their central supports in $\EHq\1{**}$
are mutually disjoint. Hence, the projections $P(\phi)$ with
$\norm{\phi}=1$, and the single projection $P(0)$ (for all the
non-\pco\ \rep s) precisely label the quasi-equivalence classes of
coherent \rep s.

{}From (3) one readily concludes that any \rep\ space $\R$ can be
split into a direct sum $\R=\R_\phi\oplus\R_\phi\1\perp$, where
$\R_\phi$ is the cyclic subspace containing $\P(\phi)\R$. Then the
\rep\ restricted to the first summand is a direct multiple of
$\pi_\phi$ with multiplicity $\dim P(\phi)\R$. The decomposition
into a Fock and a non-Fock sector (of which Lemma 3 is a special
case) is obtained for $\phi=0$. It is especially useful because the
orthogonal complement $\R_0\1\perp$ has an interesting description
\cite{\GPOTS,\WICK}:
it consists of all vectors with an ``infinite iteration history''
with respect to the operators $\as(f)$, i.e.\ it is the intersection
over $n\in\Nl$ of the closed subspaces generated by all
vectors of the form $\as(f_1)\cdots\as(f_n)\Phi$ with $f_1,\ldots,
f_n\in\R$, and $\Phi\in\R$. This decomposition can be viewed as an
analogue of the ``Wold decompositon'' of a contraction operator in
Hilbert space \cite\Foias.

\beginsection 4. Conjecture C and the Fock \rep

We begin by making precise the Conjecture C mentioned in the
introduction. We present it here, not because we are completely
convinced of its truth, but because we believe that it presents an
excellent target for future research.

\iproclaim 10 Conjecture C.
Let $-1<q<1$, and let $\H$ be a Hilbert space with $\dim\H=d<\infty$.
Let $\EHq$, and $\EH0$ denote the universal algebras introduced in
Proposition 1, and denote the respective generators by
$\as(f)\in\EHq$, and $\vs(f)\in\EH0$.
Let
$$\rho=\left({\textstyle\sum_{i=1}\1d
                        \as(e_i)\a(e_i)}\right)\1{{1\over2}}
         \in\EHq $$
for some (or any) orthonormal basis $e_1,\ldots,e_d\in\H$.
\hfill\break
Then there is a C*-isomorphism $\isom:\EH0\to\EHq$ such that
$$ \as(f)= \rho\ \isom\bigl(\vs(f)\bigr)
\quad.$$
Moreover, $0$ is an isolated point in the spectrum of $\rho$, and
the eigenprojection corresponding to this eigenvalue is
$\isom\bigl(\idty-\sum_{i=1}\1d\vs(e_i)v(e_i)\bigr)$.
\eproclaim

Note that this Conjecture can only be formulated for finite $d$,
since the sum defining $\rho$ cannot converge in norm (even though
it converges strongly in every \rep). For this reason the universal
C*-algebras for the case of infinitely many generators have to be
treated separately. For $q=0, d=\infty$ it is well known that
$\EH0\cong\OH0$ is simple, whereas it has an ideal
isomorphic to the compact operators for $d<\infty$. Analogous
phenomena occur for $q\neq0$, at least in the Fock \rep\
\cite\Burkhard. There are a number of interesting equivalent
reformulations of the Conjecture. The following one is the form in
which this Conjecture was proven \cite\QCR\ for all finite $d$, and
the restricted range $\abs q<\sqrt2-1$.

\iproclaim 11 Proposition.
Let $-1<q<1$, and let $\H$ be a Hilbert space with
$\dim\H=d<\infty$, and let $e_1,\ldots,e_d\in\H$ be an orthonormal
basis.
Then Conjecture C is equivalent to the conjunction of the following
two statements:
\item{(A)}
In the C*-algebra $\M_d(\EHq)$ of $d\times d$-matrices with entries
in $\EHq$, the matrix $X_{ij}=\a(e_i)\as(e_j)$ is strictly positive.
\item{(B)}
Let $\R$ be a Hilbert space, and let $v_i, i=1,\ldots,d$ be bounded
operators on $\R$ satisfying the relations
$v_iv_j\1*=\delta_{ij}\idty$. Then there is a unique positive
semidefinite bounded operator $\rho$ on $\R$ such that
$\as(e_i)=\rho v_i\1*$ satisfies relations (1), and such that
$\idty-\sum v_i\1*v_i$ projects onto the kernel of $\rho$.
Moreover, this unique $\rho$ necessarily lies in the C*-algebra
generated by the operators $v_i$.
\eproclaim

\proof:%
\def\As{A\1\dagger}%
\def\Vs{V\1\dagger}%
Assume (A). Consider in the universal \rep\
$\pi:\EHq\to\B(\R)$ the operators
$$\eqalign{
       \As&:\H\otimes\R     \to\R \cr
       \As&: f\otimes \psi  \mapsto \pi\bigl(\as(f)\bigr) \psi \cr
        A &:\R              \to  \H\otimes\R    \cr
        A &: \psi           \mapsto \textstyle
                     \sum_{i=1}\1d  e_i\otimes\a(e_i)\psi
\quad.}$$
Then $\rho\12=\As A$, and $X=A\As$.
The polar decomposition of $\As$ takes the form
$\As=\rho \Vs=\Vs X\1{1/2}$. Since $X>0$, the components $v_i$
of $\Vs:f\otimes\psi\mapsto=\sum_i \bra e_i,f> v_i\psi$
are in the C*-algebra $\pi(\EHq)$, and $V\Vs=\idty_{\H\otimes\R}$.
The latter relation translates into $v_iv_j\1*=\delta_{ij}\idty$.
With $\vs(f):=\sum_i \bra e_i,f> v_i\1*$, these are the $q$-relations
for $v$ with $q=0$. Hence by the universal property there is a
homomorphism $\isom:\EH0\to\EHq$ with the required property. Since
$\EH0$ has only one proper two-sided ideal, and this ideal is clearly
not annihilated by $\isom$ (consider the Fock \rep\ of $\EHq$),
$\isom$ is injective. It remains to be shown that $\isom$ is onto.
This readily follows from condition (B).

Conversely, assume that Conjecture C holds. Then in the universal
\rep\ of $\EHq$ the isomorphism $\isom$ provides a polar
decomposition of the operator $\As$. Since the polar isometry in
this case is an isometry, we must have that $A\As$ has no kernel.
If the spectrum of $A\As$ in $\M_d(\EHq)$ had an accumulation point
at zero, zero would also have to be an eigenvalue by compactness of
the state space, and the universality of the \rep. Hence
the spectrum must be bounded away from zero (A). Suppose that
$v_i$ and $\rho$ are as in (B). Then by the universality of $\EHq$
there is a unique *-\rep\ $\Phi:\EHq\to\B(\H)$ such that
$\Phi(\as(e_i))=\rho v_i\1*$. The polar decomposition of $\As$ in this
\rep\ is given by $\rho$ and $v_i$. On the other hand, by the
isomorphism with $\EH0$ we know that
$\rho=\Phi(\sum_i\as(e_i)\a(e_i))$ is in the C*-algebra generated by
the $v_i\1*=\Phi(\vs(e_i))$.
\QED

Some consequences of Conjecture C would be the following: (1) The
Fock \rep\ of $\EHq$ is faithful for all $q$. (2) $\EHq$ has only
one proper ideal, isomorphic to the compact operators. (3) the
resulting quotient is isomorphic to the Cuntz algebra
$\OH0$. Statement (3) may be extended to a version of
Conjecture C on the level of the $q$-Cuntz algebras $\OHq$. Since
$\OHq$ can be obtained as a quotient of {\it any} other \rep\ of
$\EHq$, we can utilize specific information about the Fock
\rep\ in approaching this problem.

Dykema and Nica \cite\Dykema, building on results of Zagier
\cite\Zagier, were able to verify parts of Conjecture C in the Fock
\rep $\pi_0$. For example, they verified statement (A) of Proposition
11 for that \rep, by unitary implementation of a homomorphism
$$ \isom_0:\pi_0(\EH0)\to\pi_0(\EHq)
\quad\eqno(16)$$
satisfying the properties required in Conjecture C, except
surjectivity. Their results imply a lower bound
$$\eqalign{
  \pi_0(X)   &\geq {1-q\over1-\absq}\ \epsilon(\absq)\ \idty
              >0  \cr
  \epsilon(s)&= \prod_{k=1}\1\infty{1-s\1k\over1+s\1k}
              = \sum_{k=-\infty}\1\infty (-1)\1k\, s\1{k\12}
}\eqno(17)$$
for $X\in\M_d(\EHq)$, $X_{ij}=\a(e_i)\as(e_j)$. (We remind the reader
of the difference in normalization between \cite\Dykema, and this
paper). Moreover, they showed surjectivity of $\isom_0$ for
$$\eqalign{
      q\12 &<\epsilon(\absq) \cr
\hbox{i.e.}\qquad
    \absq &<\approx 0.44
}\quad.\eqno(18)$$
We can immediately translate these results into a partial
verification of Conjecture C on the Cuntz algebra level:

\iproclaim 12 Theorem.
Let $-1<q<1$, and let $\H$ be a Hilbert space with $\dim\H=d<\infty$.
Let $\OHq$, and $\OH0\equiv\Cuntzal_d$ denote the $q$-Cuntz algebra,
and the Cuntz algebra, as in Definition 8,
and denote by $\pi_1:\EHq\to\OHq$ (or $\EH0\to\OH0$)
the respective quotient maps.
Let $\rho\in\EHq$ be as in Conjecture C.
\item{(1)}
Then there is a (not necessarily surjective)
C*-homomorphism $\isOm:\OH0\to\OHq$ such that
$\pi_1(\as(f))= \pi_1(\rho)\ \isOm\pi_1\bigl(\vs(f)\bigr)$.
\item{(2)}
$$ \pi_1(\rho\12)\geq {1-q\over1-\absq} \epsilon(\absq)\idty
\quad.$$
\item{(3)}
Let $\omega_\phi$ be a \pco\ state on $\OHq$.
Then $\omega_\phi \circ \isOm$ is the \pco\ state on $\OH0$
associated with $\phi$.
\item{(4)}
When $q\12<\epsilon(\absq)$, $\isOm$ is onto,
and hence an isomorphism.
\eproclaim

\proof:
The eigenprojection onto the kernel of $\pi_0(\rho)$ is
$P_0 \equiv\isom_0(\idty-\pi_0(\sum v_i\1*v_i))
      \in\pi_0(\EHq)$.
Consider a \pco\ \rep\ $\pi_\phi$ of $\pi_0(\EHq)$. Since \pco\
states are pure, this \rep\ is irreducible. On the other hand,
$\pi_\phi(P_0)$ projects onto the set of Fock vectors in that \rep.
The invariant subspace generated from a Fock vector is a copy of
Fock space, on which the projection $P(\phi)$, as in Proposition 9,
vanishes. On the other hand, $P(\phi)\neq0$, so the Fock sector
cannot be the whole space, and must be zero by irreducibility. Hence
$\pi_\phi(P_0)=0$, and $\pi_\phi(\rho)>0$. The bound (2) then
follows from equation (17). Moreover, under the map
$$\def\Ra#1{\ {\buildrel{#1}\over\longrightarrow\ }}
     \EH0\Ra{\pi_0}   \pi_0(\EH0)
         \Ra{\isom_0} \pi_0(\EHq)
         \Ra{\pi_\phi} \pi_\phi(\EHq)\cong\OHq
$$
$\idty-\sum_i\vs(e_i)v(e_i)$ becomes $\pi_\phi(P_0)=0$. Hence it
lifts to the quotient as $\isOm:\OH0\to\EH0$. The properties (1),
(4) of this map are readily verified from those proven for the Fock
\rep.

To see (3), recall from Lemma 3 that the eigenvalue equation
$\pi_\phi\bigl(\a(\phi)-1\bigr)\omega_\phi=0$ can only be satisfied
when we also have $\pi_\phi\bigl(\as(\phi)-1\bigr)\omega_\phi=0$.
Hence with a basis $e_1=\phi,e_2,\ldots,e_d\in\H$ we get
$$ \pi_\phi(\rho\12)\Omega_\phi
     = \sum_i \pi_\phi\bigl(\as(e_i)\a(e_i)\bigr)\Omega_\phi
     = \pi_\phi\bigl(\as(e_1)\a(e_1)\bigr)\Omega_\phi
     = \Omega_\phi
\quad.$$
Since $\pi_\phi(\rho)>0$, this entails
$$\eqalign{
 \pi_\phi\isOm\bigl(v(f)\bigr)\Omega_\phi
    =\pi_\phi(\rho)\1{-1}\ \pi_\phi(\a(f))\Omega_\phi
    = \bra f,\phi>\ \pi_\phi(\rho)\1{-1} \Omega_\phi
    = \bra f,\phi> \Omega_\phi
\quad.}$$
Therefore, for $X\in\EH0$,
$\omega_\phi\bigl(\isOm(Xv(f))\bigr)
      = \bra f,\phi>\omega_\phi\bigl(\isOm(X)\bigr)$,
which proves (3).
\QED

\beginsection 5. The boundary points $q=\pm1$

Apart from Conjecture C and its special cases, an interesting
problem concerning the $q$-relations (1) is to show that they define
a continuous field of C*-algebras $\EHq$ in the parameter $q$ in the
sense of Dixmier \cite\Dix. If Conjecture C holds, i.e.\ an
isomorphism $\isom_q:\EHq\to\EH0$ exists, this problem amounts, for
$q\neq\pm1$, to the question whether the element
$\isom_q\1{-1}(\rho)\in\EH0$ depends continuously on $q$. (For
$\absq<\sqrt2-1$, this continuity is easily verified from the
argument in \cite\QCR). The interesting questions arise at the
boundaries $q=\pm1$.

The role of the coherent states is that of a {\it continuous field of
states} in the following sense: for any polynomial $X$ in the
variables $\as(f)$, $\a(g)$, and $q$ (with $f,g\in\H$), and for every
fixed $\phi\in\H$, the coherent expectation $\omega_\phi(X)$ is a
continuous function of $q$. The continuity of the field $q\mapsto\EHq$
is related to the existence of sufficiently many such continuous
fields of states.

As a first step towards understanding the continuity at $q=\pm1$, we
compute the algebras $\EH{\pm1}$, and their coherent states. Recall
that for $q=1$ we have imposed, in Proposition 1,  the bound
$\sum_i\as(e_i)\a(e_i)\leq\idty$ for any family of orthogonal
vectors.

\iproclaim 13 Proposition. Let $\H$ be a Hilbert space. Then $\EH{+1}$
is isomorphic to the algebra of weakly continuous functions on the
unit ball of $\H$. A state on this algebra is coherent if and only
if it is pure.
\eproclaim

\proof:
We have to show first that $\EH{+1}$ is abelian. Clearly,
$\bracks{\as(f),\a(g)}=0$ for all $f,g$. In particular, each $\as(f)$
is a bounded normal operator. By Fuglede's Theorem
\cite{\Fuglede,\Rudin}, $\as(f)$ and $\as(g)$ also commute, and
$\EH{1}$ is abelian. A pure state $\omega$ must be multiplicative,
and is hence determined by its value $\omega(\as(f))$ on the
generators. Since $\as$ is linear, and since
$\as(f)\a(f)\leq\norm{f}\12\idty$, this expression must be a bounded
linear functional, and hence of the form
$\omega(\as(f))=\bra\phi,f>$ with $\phi\in\H, \norm{\phi}\leq1$.
Hence $\omega=\omega_\phi$ is coherent, and any coherent state is
obtained in this way. Note that, for all polynomials $X$ in the
generators $\as(f),\a(f)$, the function $\phi\mapsto\omega_\phi(X)$
is weakly continuous on the unit ball. On the other hand, the
algebra of such polynomials is dense in the algebra of all weakly
continuous functions by the Stone-Weierstra\ss\ Theorem.
\QED

Note that by this Proposition the set of coherent states is faithful
at $q=+1$. This suggests that they may be useful for proving the
continuity at $q=1$, provided one can show that collection of coherent
\rep s is also faithful for $q<1$. In the following Proposition we
see that faithfulness does not hold at the other limit point $q=-1$,
where the relations become
$$ \a(f)\as(g)+\as(g)\a(f)=2\,\bra f,g>\idty
\quad.\eqno(21)$$
The Proposition is based on well-known results in the theory of
Clifford algebras \cite{\Bourbaki,\Shale,\Lawson}, which arise from
these relations either by taking $f,g$ to be in a real Hilbert
space, and setting $\as(f)=\a(f)$. In even dimensions this is
equivalent to taking (21), and adding the relation that the $\as(f)$
anti-commute with each other. The algebra arising in this way is the
Fock \rep\ of (21), and we will refer to it as the CAR-algebra
\cite\BraRo. The point of the following Proposition is that no
anti-commutation relation is added, but that such a relation
automatically holds in every \irrep.

\iproclaim 14 Proposition. Let $\H$ be a Hilbert space, and consider
the C*-algebra $\EH{-1}$ as defined in Proposition 1. Then
\item{(1)}
The elements
$$ \as(f)\as(g)+\as(g)\as(f) =2\thetah(f,g)
\quad,$$
for $f,g\in\H$ generate the center of $\EHq$.
\item{(2)}
The center of $\EH{-1}$ is isomorphic to $\C(\St)$, where $\St$ is the
set of all symmetric bilinear forms $\theta:\H\times\H\to\Cx$ such
that
$$  \abs{\theta(f,g)}\leq \norm{f}\,\norm{g}    $$
for all $f,g\in\H$, with the coarsest topology making the functions
$\theta\mapsto\theta(f,g)$ continuous.
\item{(3)}
Let $\theta$ be a symmetric bilinear form satisfying the above
bound, and let $\N(\theta)$ denote the real subspace of vectors
$f\in\H$ such that $\theta(f,f)=\norm{f}\12$. Let $r(\theta)$ denote
the dimension of the complement of $\N(\theta)$ in $\H$, taken as a
real Hilbert space.
Let $\EH{-1,\theta}$ denote the quotient of $\EH{-1}$ by the
relations $\thetah(f,g)=\theta(f,g)\idty$. Then
\itemitem{(a)} If $r(\theta)$ is finite and even, $\EH{-1,\theta}$
is isomorphic to the algebra of $2\1{r(\theta)/2}$-dimensional
matrices.
\itemitem{(b)} If $r(\theta)$ is finite and odd, $\EH{-1,\theta}$
is isomorphic to the direct sum of two copies of the algebra of
$2\1{(r(\theta)-1)/2}$-dimensional matrices.
\itemitem{(b)} If $r(\theta)$ is infinite, $\EH{-1,\theta}$ is
isomorphic to the CAR-algebra on an infinite dimensional Hilbert
space.
\eproclaim

\proof:%
\def\Stt{\widetilde\St}%
\def\sh{\widehat s}%
\def\thetat{\widetilde\theta}
(1) By an elementary computation one verifies that
$\thetah(f,g)$ commutes with all $\a(h)$. Hence
$\thetah(f,g)$ is normal, and by Fuglede's Theorem \cite\Rudin\ it
must also commute with $\as(h)$. Hence $\thetah(f,g)$ is in the
center for all $f,g\in\H$. Let $\C(\St)\subset\EH{-1}$ denote the
C*-algebra generated by the $\thetah(f,g)$. Its spectrum space $\St$
is the set of those symmetric bilinear forms $\theta$, which may
arise in an \irrep\ of $\EH{-1}$, i.e.\ those $\theta$ for which the
relations
$$\eqalign{
       \a(f)\as(g)+\as(g)\a(f)&=2\,\bra f,g>\idty \cr
    \as(f)\as(g)+\as(g)\as(f) &=2\theta(f,g)\idty
\quad}\eqno(*)$$
have a solution by a bounded linear operator $\as:\H\to\B(\R)$ for
some Hilbert space $\R$. The rest of the proof depends on the analysis
of this set of relations.

The unique feature of the relations (1) in at $q=-1$ is that the
symmetry with respect to exchange of $\a$ and $\as$ (up to questions
of linearity/antilinearity). Therefore we will consider $\H$ now as a
real Hilbert space of dimension $\dim_\Rl(\H)=2\dim_\Cx(\H)$, and
introduce the hermitian generators
$$ s(f):={1\over2}\bigl(\as(f)+\a(f)\bigr)
\quad,$$
which are real linear in $f\in\H$. From $s(f)$ and the complex
structure on $\H$ we can recover the original generators by the
formula $\as(f)=s(f)-i\,s(if)$. In terms of the new generators we get
the relations
$$ s(f)s(g)+s(g)s(f)= \Re\Set\Big{ \bra f,g>+ \theta(f,g)}
                    =: 2 \Theta(f,g)\idty
\quad.\eqno(**)$$
Clearly, $\Theta$ is a symmetric, real valued form on the real
Hilbert space $\H$. Since $\Theta(f,f)=s(f)\12\geq0$, the positivity
of $\Theta$ is necessary for the existence of a \rep, and hence for
$\theta\in\St$.

We claim that the positivity of $\Theta$ is equivalent to the
inequality in item (2) of the Proposition. Clearly,
$\abs{\theta(f,f)}\leq\norm{f}\12$ is sufficient for
$\Theta(f,f)\geq0$. Conversely, assume
$2\Theta(f,f)=\norm{f}\12+\Re\theta(f,f)\geq0$. Substituting $f\mapsto
if$ in this inequality we get that
$\abs{\Re\theta(f,f)}\leq\norm{f}\12$. Hence by the Schwarz inequality
in the real Hilbert space $\H$ we have
$\abs{\Re\theta(f,g)}\leq\norm{f}\norm{g}$, and the result follows by
replacing $f$ in this inequality by a complex multiple of $f$.
It is also easy to see that the rank of $\Theta$ is equal to
$r(\theta)$, as defined in item (3).

In order to prove the characterization of $\St$, the joint spectrum of
the central elements $\thetah(f,g)$, it remains to be shown that for
every $\Theta\geq0$ there is some \rep\ of $(*)$. We will
simultaneously prove (3) by constructing all such \rep s (assuming
$\Theta\geq0$), and showing that they have the form given in (3) with
$r(\theta)=\rank\Theta$.

We can find an orthonormal basis $\set{e_i}\subset\H$ such that
$\Theta(e_i,e_j)=\Theta_i\delta_{ij}$. The generators $s(e_i)$ with
$\Theta_i=0$ have to be zero, and the remaining ones can be multiplied
by $\Theta_i\1{-1/2}$, so that $(**)$ becomes equivalent to the
relations
$$    s_is_j+s_js_i= 2\, \delta_{ij}\, \idty
\quad,\eqno(***)$$
where the $s_i=s_i\1*$, and $i=1,\ldots,\rank\Theta$.
Hence the isomorphism type of $\EH{-1,\theta}$ depends only on
$\rank\Theta$. One readily verifies that $\rank\Theta=r(\theta)$, as
defined in (3).

For finite $r(\theta)$, (3) follows from the standard results of
the \rep\ theory of Clifford algebras (see e.g.\ Theorems 2 and 3 in
\S9,No.4 of \cite\Bourbaki). The CAR-algebras are a special case of
these arguments with $\theta=0$. For infinitely many generators $s_i$
we therefore get an inductive limit which we can take along the
simple algebras of even numbers of generators \cite\Shale, and which
is identical with the inductive limit defining the CAR-algebra over an
infinite dimensional Hilbert space. We note that in the case (3b) the
center is generated by the odd element
$$ \sh=s_1\cdots s_{r(\theta)}
\quad,$$
which is unitary and satisfies $\sh\12=\pm\idty$, depending on
$r(\theta)$ modulo $4$ \cite\Lawson. In any case, $\sh$ has two
eigenvalues $\pm1$ or $\pm i$, which label the two \irrep s with
given $\theta$, and are exchanged by the parity automorphism defined
by $\as(f)\mapsto-\as(f)$.

We have shown that the algebra generated by the elements $\thetah$ is
isomorphic to $\C(\St)$ with $\St$ as described in item (2). It
remains to be shown that this algebra coincides with the center of
$\EH{-1}$. Let the center of $\EH{-1}$ be $\C(\Stt)$ for some compact
space $\Stt$. Since $\C(\St)$ is a subalgebra, we have a canonical
continuous surjection $p:\Stt\to\St$. Whenever $r(\theta)$ is even the
relations $(*)$ have only one \irrep, which implies that
$p\1{-1}(\set\theta)$ is a single point. Otherwise,
$p\1{-1}(\set\theta)$ may consist of at most two points, corresponding
to the two \irrep s of $(*)$. The parity automorphism induces a
homeomorphism $F:\Stt\to\Stt$ which leaves all points with even
or infinite $r(\theta)$ fixed. Whenever $r(\theta)$ is odd, and
$p\1{-1}(\set\theta)$ consists of two points, these two points are
exchanged by $F$.

Since $\H$ has even or infinite real dimension, $r(\theta)$ is
odd or infinite for a dense subset of $\theta\in\St$. Now consider
some $\theta$ with odd $r(\theta)$, and let $\theta_\alpha\in\St$
be a net with $\theta_\alpha\to\theta$, and $r(\theta_\alpha)$ even
for all $\alpha$. Let $\thetat_\alpha$ be the net in $\Stt$ uniquely
defined by $p\bigl(\thetat_\alpha\bigr)=\theta_\alpha$. Since $F$ is
continuous, and $F\thetat_\alpha=\thetat_\alpha$ any cluster point
$\thetat$ of this net must also be fixed under $F$, and since $p$ is
continuous, we must have $p(\thetat)=\theta$. But the only way
$\thetat\in p\1{-1}(\set\theta)$  can be fixed by $F$ is that
$p\1{-1}(\set\theta)$ is a single point.
It follows that $p:\Stt\to\St$ is a bijection, and the center of
$\EH{-1}$ coincides with the algebra $\C(\St)$ generated by the
$\thetah(f,g)$.
\QED

It is clear that the coherent \rep s of $\EH{-1}$ are precisely
those for which
$$ \theta(f,g)=\bra\phi,f>\,\bra\phi,g>
\quad\eqno(22)$$
is a rank one operator. The set $\N(\theta)$ of vectors $f$ with
$\norm{f}\12=\theta(f,f)$ is either null, when $\norm{\phi}<1$, or is
the the one-dimensional real subspace spanned
by $\phi$ when $\norm{\phi}=1$. Hence for the \pco\ states on
$\EH{-1}$ with $\dim\H<\infty$, $r(\theta)$ is odd.

When $\dim\H=1$, all symmetric bilinear forms on $\H$ are of the form
(21). Hence in this case the set of coherent states provides an
everywhere faithful family of continuous fields of states.
Accordingly, $q\mapsto\EHq$ is a continuous field of C*-algebras
\cite\Nica. For $\dim\H>1$ an interesting problem arises here: since
the rank one bilinear forms are a low dimensional subset of $\St$ it
is clear that many \irrep s of $\EH{-1}$ are not coherent \rep s. Is
it possible to embed states on such non-coherent \rep s of $\EH{-1}$
into a continuous field of states for the field $q\mapsto\EHq$?

\refskip=3.7pt \let\REF\doref

\ACKNOW
This paper grew out of a collaboration \cite\QCR\ with Lothar Schmitt
(Osnabr\"uck). It has benefited from conversations with R.~Speicher,
B.~K\"ummerer, A.~Nica, and K.~Dykema, whom we also thank for making
their work available to us. R.F.W.\ acknowledges financial support
from the DFG (Bonn) through a scholarship and a travel grant.

\REF Ask \Askey \Gref
    R. Askey
    "Continuous $q$-Hermite polynomials when $q>1$" pp. 151-158
    \hfill\break \inPr
    D. Stanton
    "$q$-series and partitions"
    Springer-Verlag, New York 1989

\REF Bae \Baez \Jref
    J.C. Baez
    "$R$-commutative geometry and quantization of Poisson algebras"
    \hfill\break
    Adv.Math. @95(1992) 61-91

\REF Ber \Bergmann \Jref
    G. Bergmann
    "The diamond lemma for ring theory"
    Adv.Math. @29(1978) 178-218

\REF Bie \Bieden \Jref
    L.C. Biedenharn
    "The quantum group $SU_q(2)$ and a $q$-analogue of the boson
    operators"
    J.Phys.A. @22(1989) L873-L878

\REF Bou \Bourbaki \Bref
     N. Bourbaki
     "\'El\'ements de math\'ematique, Alg\`ebre, Chapitre 9"
     Hermann, Paris 1959

\def\Bozejko{Bo{\accent 95 z}ejko}
\REF BS1   \BoSpeia  \Jref
    M. \Bozejko, R. Speicher
    "An example of a generalized Brownian motion"
    \hfill\break
    Commun.Math.Phys. @137(1991) 519-531
    \more; \hfill\break
    part II
    \inPr L. Accardi "Quantum Probability \& related topics VII"
    World Scientific, Singapore 1992

\REF  BEGJ \JOa \Jref
    O. Bratteli, D.E. Evans, F.M. Goodman, P.E.T. J\o rgensen
     "A dichotomy for derivations on ${\cal O}_n$"
    Publ. RIMS, Kyoto Univ. @ 22(1986) 103--117

\REF BR \BraRo      \Bref
    O. Bratteli, D.W. Robinson
     "Operator algebras and quantum statistical mechanics"
     volume II, Springer Verlag, Berlin, Heidelberg, New York
    1981

\REF Cun \Cuntz  \Jref
    J. Cuntz
    "Simple C*-algebras generated by isometries"
    \hfill\break
    Commun.Math.Phys. @57(1977) 173-185

\REF Dix \Dix \Bref
     J. Dixmier
     "C*-algebras"
     North-Holland, Amsterdam 1982

\REF DS \Dunford \Bref
     N. Dunford, J.T. Schwarz
     "Linear Operators, I: General Theory"
     Interscience, New York 1958

\REF DN \Dykema  \Gref
    K. Dykema, A. Nica
    "On the Fock representation of the $q$-commutation relations"
    Preprint Berkeley 1992

\REF Fiv \Fivel  \Jref
     D.I. Fivel
     "Interpolation between Fermi and Bose statistics using
      generalized commutators"
      Phys.Rev.Lett. @65(1990) 3361-3364\more
      .\ Erratum: \Jn
      Phys.Rev.Lett. @69(1992) 2020

\REF Fug \Fuglede \Jref
    B. Fuglede
    "A commutativity problem for normal operators"
    Proc.Nat.Acad.Sci. USA @36(1950) 35-40

\REF Gre \Green \Jref
     O.W. Greenberg
     "Particles with small violations of Fermi or Bose statistics"
     \hfill\break
     Phys.Rev.D. @43(1991) 4111-4120

\REF GM \Greenloc \Jref
    O.W. Greenberg, R.N. Mohapatra
    "Difficulties with a local quantum field theory of possible
    violation of the Pauli principle"
    Phys.Rev.Lett. @62(1989) 712-714

\REF HR \Rieckers \Jref
     R. Honegger, A. Rieckers
     "The general form of non-Fock coherent Boson states"
     Publ.RIMS, Kyoto Univ. @26(1990) 945-961

\REF JSW1 \QCR \Gref
     P.E.T. J\o rgensen, L.M. Schmitt, R.F. Werner
     "$q$-canonical commutation relations
       and stability of the Cuntz algebra"
     Preprint Osnabr\"uck 1992;
     to appear in {\it Pacific.J.Math.}

\REF JSW2 \GPOTS \Gref
     P.E.T. J\o rgensen, L.M. Schmitt, R.F. Werner
    "$q$-relations and stability of C*-isomorphism classes"
    To appear in the proceedings of the {\it Great Plains operator
    theory symposium}, Iowa 1992

\REF JSW3 \WICK \Gref
     P.E.T. J\o rgensen, L.M. Schmitt, R.F. Werner
    "Positive representations of general Wick ordering commutation
     relations"
     In preparation

\REF KaS \Katriel \Jref
     J. Katriel, A.I. Solomon
     "A $q$-analogue of the Campbell-Baker-Hausdorff expansion"
     J.Phys.A. @24(1991) L1139-L1142

\REF KlS \Klauder \Bref
    J.R. Klauder, E.C.G. Sudarshan
    "Fundamentals of Quantum Optics"
    Benjamin, New York 1968

\REF K\"uS1 \KuSpei \Jref
    B. K\"ummerer, R. Speicher
    "Stochastic integration on the Cuntz algebra ${\cal O}_\infty$"
    J.Funct.Anal. @103(1992) 372-408

\REF K\"uS2  \Burkhard
     \authors{B. K\"ummerer, R. Speicher, M. \Bozejko},
      unpublished results

\REF LM \Lawson \Bref
    H.B. Lawson{, Jr.}, M.-L. Michelson
    "Spin geometry"
    Princeton Univ. Press, Princeton 1989

\REF LS \LiSheng \Jref
    Y-Q. Li, Z-M. Sheng
    "A deformation of quantum mechanics"
    J.Phys.A @25(1992) 6779-6788

\REF Mcf \Farlane \Jref
     A.J. Macfarlane
     "On $q$-analogues of the quantum harmonic oscillator and the
     quantum group $SU(2)_q$"
     J.Phys.A. @22(1989) 4581-4588

\REF NN \Nica  \Gref
     G. Nagy, A. Nica
     "On the `quantum disk' and a `non-commutative circle' "
     Preprint, Berkeley 1992.
     To appear in the proceedings of the {\it Great Plains operator
     theory symposium}, Iowa 1992

\REF  S-NF \Foias \Bref
    B. Sz.-Nagy, C. Foia\c s
    "Harmonic analysis of operators on Hilbert space"
    North-Holland Publ., Amsterdam 1970

\REF Rud \Rudin \Bref
    W. Rudin
    "Functional analysis"
    McGraw-Hill, 1975

\REF Sak \Sakai \Bref
    S. Sakai
    "C*-algebras and W*-algebras"
    Springer-Verlag, Ergebnisse 60, New York 1971

\REF SS \Shale \Jref
    D. Shale, W. Stinespring
    "States of the Clifford algebra"
    Ann.Math. @80(1964)365-381

\REF Sp1  \Speia  \Jref
    R. Speicher
    "A new example of `Independence' and `white noise' "
    \hfill\break
    Probab.Th.Rel.Fields @84(1990)141-159

\REF Wer \QTD \Gref
    R.F. Werner
    "The Free Quon Gas Suffers Gibbs' Paradox"
    Preprint, Osnabr\"uck 1993

\REF Wor \Worono  \Jref
    S.L. Woronowicz
    "Twisted SU(2) group. An example of non-commutative differential
    calculus"
    Publ.RIMS, Kyoto Univ. @23(1987)117-181

\REF Zag \Zagier \Jref
    D. Zagier
    "Realizability of a model in infinite statistics"
    \hfill\break
    Commun.Math.Phys. @147(1992) 199-210

\bye